\documentclass[12pt,intlim,righttag]{article}
\usepackage{amsmath,amsthm}
\usepackage{amssymb}
\usepackage{amsfonts}

\newcommand{\la}{\langle}
\newcommand{\ra}{\rangle}

\newcommand\beq{\begin{equation}}
\newcommand\eeq{\end{equation}}

\theoremstyle{Theorem}

\theoremstyle{corollary}

\theoremstyle{remark}

\theoremstyle{definition}

\begin{document}
\title{Remark on Right Continuous Exponential Martingales}

\author{B. Chikvinidze }

\date{~}
\maketitle

\begin{center}

  Institute of Cybernetics of Georgian Technical University,
  
  \
  \\       
  Georgian American University, Business school, 8 M. Aleksidze Srt., \\
          Tbilisi 0160, Georgia 

\          
\\
   E-mail: beso.chiqvinidze@gmail.com \\

\begin{abstract}


{\bf Abstract}     Using $\la M^c \ra $, jump measure $\mu$ and its compensator $\nu$ we characterize the event where the stochastic exponential $\mathcal{E}(M)$ equals to zero.

\end{abstract}

\bigskip


\


\end{center}

\
\\
{\bf 1. Introduction. }Let us introduce a basic probability space $\big( \Omega , \mathcal{F}, P \big)$ and a right continuous filtration 
$(\mathcal{F}_t)_{0\leq t < \infty}$ satisfying usual conditions. Let $\mathcal{F}_\infty $ be the smallest $\sigma-$Algebra containing all $\mathcal{F}_t $ for $t<\infty$ and let $M = (M_t)_{t \geq 0}$ be a local martingale on the stochastic interval $[[0;T[[$, where $T$ is a stopping time. Denote 
by $\bigtriangleup M_t = M_t - M_{t-}$ jumps of $M$ and
by $\mathcal E(M)$ the stochastic exponential of the local martingale $M$:

$$
\mathcal E_{t}(M)= \exp\big\{ M_t-\frac{1}{2}\la M^c\ra_t \big\} \prod_{0<s \leq t}(1+\bigtriangleup M_s)e^{-\bigtriangleup M_s},
$$
where $M^c$ denotes a continuous local martingale part of $M$. Notice, that $M=M^c + M^d$ where $M^d$ is a purely discontinuous local martingale part of $M$, which means that $M^d$ is orthogonal to any continuous local martingale. With this we known that $M^d_t=\int_0^t \int_{-1}^{\infty} x d(\mu -\nu)$, where $\mu(\omega,t,x)$ is the jump measure of $M$ and 
$\nu(\omega,t,x)$ is it's compensator. 

\
Through this paper we will integrate with respect to $\mu $ over the set $(-1;1)\setminus\{0\}$ and we will write it as $\int_0^T \int_{-1}^{1}\cdot \; d\mu $.

\
It is well known that $\mathcal E_{t}(M) = 1 + \int^t_0 \mathcal E_{s-}(M) d M_s$, so it is clear that for local martingale $M$ the associated stochastic exponential $\mathcal E(M)$ is a local martingale. Throughout of this paper we assume that 
$\bigtriangleup M_t \geq -1$ which implies that $\mathcal E(M)$ is a non-negative local martingale and therefore a supermartingale. In case when $\mathcal{E}(M)$ is a uniformly integrable martingale on $[[0;T]]$, we can define using $\mathcal{E}(M)$ and the Radon-Nikodym derivative a new probability measure: $dQ=\mathcal{E}_T(M)dP$. It is clear that $Q<<P$ and if $P\{\mathcal{E}_T(M)>0\}=1$, then $P$ and $Q$ will be equivalent probability measures ($P\sim Q$). To know whether 
$P\sim Q$ or not, we must study the set $\{\mathcal{E}_T(M)=0\}$. In case when $M=M^c$ it was shown by Kazamaki \cite{2} in 1994 that $\{\mathcal{E}_T(M^c)=0\}=\{\la M^c \ra_T =\infty \}$. For general $M$, in 1978 it was proved by J. Jacod \cite{1} that 
$$\{ \mathcal{E}_{\infty}(M)>0 \}=\Big\{ \la M^c \ra _{\infty} + \int _0^{\infty} \int_{-1}^{\infty} \frac{x^2}{1+|x|} d\nu  + \int_0^{\infty} \frac{1}{\mathcal{E}_{s-}(M)}dB_s  < \infty  \Big\} $$  
where $B_s$ is the predictable, non-decreasing process from the Doob-Meyer decomposition of $\mathcal{E}(M)$. In 2019 M. Larsson and J. Ruf \cite{3} proved the set inclusion 
$$\{\lim_{t\uparrow \tau}\mathcal{E}_{t}(M)=0\}\subset
\big\{ \lim_{t\uparrow \tau}M_t=-\infty \big\} \cup \{[M]_{\tau } =\infty \} \cup 
\{ \bigtriangleup M_t=-1, t\in [0;\tau ) \}$$ 
holds true for any predictable stopping time $\tau $. With this they proved, that if in addition $\bigtriangleup M\geq-1$ and 
$\overline{\lim}_{t\uparrow \tau}M_t<\infty $, then the reverse set inclusion also holds. 

\
\\
The aim of this paper is to characterize the set $\{\mathcal{E}_T (M)=0\}$ using $\la M^c \ra $, $\mu(\omega,t,x)$ and 
$\nu(\omega,t,x)$, for any stopping time $T$.

\newpage
\
\\
{\bf Theorem 1} \; Let $M$ be a local martingale with $\bigtriangleup M \geq -1$. Then the following set equalities hold true $P$ a. s.:  
$$(i) \;\;\; \big\{ \mathcal E_{T}(M)=0 \big\}=\Big\{ \la M^c \ra_{T} + \int_0^T \int_{-1}^1 \frac{x^2}{1+x}d\mu + \int_0^T \int_{1}^{+\infty}\frac{x^2}{1+x}d\nu = \infty \Big\};$$
\\
$(ii)$ \; If $E\frac{1}{1+\bigtriangleup M_{\sigma}}1_{\{|\bigtriangleup M_{\sigma}|\leq 1\} }<\infty $, for any $\sigma<\infty $, then
$$\big\{ \mathcal E_{T}(M)=0 \big\}=\Big\{ \la M^c \ra_{T} + \int_0^T \int_{-1}^{+\infty}\frac{x^2}{1+x}d\nu = \infty \Big\};
$$
\\
$(iii)$ \; If $E\bigtriangleup M_{\sigma}<\infty $, for any $\sigma<\infty $, then
$$\big\{ \mathcal E_{T}(M)=0 \big\}=\Big\{ \la M^c \ra_{T} + \int_0^T \int_{-1}^{+\infty}\frac{x^2}{1+x}d\mu = \infty \Big\}.$$

\
\\
{\bf Remark 1} In the contrary to the result from Jacod \cite{1}, we are not using the additional increasing process $B_t$, which is not in terms of $M$. In their result Larsson and Ruf \cite{3} used the predictable stopping time $\tau $ and they have additional restriction on $M$ to obtain the set equality. In part $(i)$ of Theorem 1 we have the set equality without any restriction on $M$ and in part $(ii)$ we have the set equality with predictable characteristics of $M$, but with integrability restriction on jumps of $M$. With this let us mention that we use any kind of stopping times $T$, while Larrson and Ruf \cite{3} used only predictable stopping times.  

\ 
\\
{\it Proof of the Theorem 1:}  If $\bigtriangleup M_s=-1$ for some $s\leq T$, then it is obvious that $\mathcal E_{T}(M^d)=0 $ and $\int_0^T \int_{-1}^1 \frac{x^2}{1+x}d\mu = \sum _{s\leq T} \frac{(\bigtriangleup M_s)^2}{1+\bigtriangleup M_s} = \infty $, so we can prove Theorem 1 when $\bigtriangleup M_s >-1$.

\
Define local martingales   
$$M^1_t = \int^t_0 \int_{-1}^1 x d(\mu - \nu);  \;\;\;\;\;   M^2_t = \int^t_0 \int_{1}^{+\infty} x d(\mu - \nu).$$ 
It is clear that $|\bigtriangleup M^1_t|\leq 1$, 
$\bigtriangleup M^2_t \geq 1$ and $M^d_t = M^1_t + M^2_t$, so we have $M = M^c + M^1 + M^2$. It is easy to check that $\mathcal{E}_T (M)=\mathcal{E}_T (M^c)\mathcal{E}_T (M^1)\mathcal{E}_T (M^2)$, so 
$$\{ \mathcal{E}_{T}(M)=0 \}=\{ \mathcal{E}_{T}(M^c)=0\}\cup \{\mathcal{E}_{T}(M^1)=0\}\cup \{\mathcal{E}_{T}(M^2)=0\}.$$ 
It is well known from Kazamaki \cite{2} that $\{ \mathcal{E}_{T}(M^c)=0\}=\{ \la M^c \ra_{T}=\infty \} $, so to prove part $(i)$ of Theorem 1 it is sufficient to show the set equalities 

\begin{equation}
\big\{ \mathcal E_{T}(M^1)=0 \big\} = \Big\{ \int_0^T \int_{-1}^1 \frac{x^2}{1+x}d\mu = \infty \Big\},
\end{equation}

\begin{equation}
\big\{ \mathcal E_{T}(M^2)=0 \big\} = \Big\{ \int_0^T \int_{1}^{+\infty} \frac{x^2}{1+x}d\nu = \infty \Big\}.
\end{equation}
 
First let us show that $\big\{ \int_0^T \int_{-1}^{+\infty} \frac{x^2}{1+x} d\mu = \infty \big\} \subset \big\{ \mathcal E_{T}(M^d)=0 \big\}$ for any local martingale $M$. An easy calculations give us:

$$\mathcal{E}^2_{T}\Big( \frac{1}{2}M^d \Big) = \exp \Big\{ M^d_{T} + \int_0^T \int_{-1}^{+\infty}\big[ 2\ln(1+\frac{x}{2})-x \big]d\mu = $$

$$\mathcal{E}_{T}(M^d)\exp \Big\{ \int_0^T \int_{-1}^{+\infty} \ln \frac{(1+\frac{x}{2})^2}{1+x}d\mu \Big\}$$ 

and from this we obtain:

$$\mathcal{E}_{T}(M^d)=\mathcal{E}^2_{T}\Big( \frac{1}{2}M^d \Big)\exp \Big\{ - \int_0^T \int_{-1}^{+\infty} \ln \Big(1+\frac{1}{4}\cdot \frac{x^2}{1+x}\Big) d\mu \Big\}.$$
\\
The supermartingale property of $\mathcal{E}\big( \frac{1}{2}M^d \big) $ implies 
$P\{ \mathcal{E}_{T}(\frac{1}{2}M^d) < \infty \} = 1 $, so we obtain that $\big\{\int_0^T \int_{-1}^{+\infty} \ln \big(1+\frac{1}{4}\cdot \frac{x^2}{1+x}\big) d\mu =\infty \big\} \subset 
\{\mathcal{E}_{T}(M^d)=0 \}$. Now the set equalities below are obvious and the first set inclusion follows from the inequality
$\ln (1+\sum_n x_n)\leq \sum_n \ln (1+x_n)$, where $x_n\geq 0$:

$$\Big\{ \int_0^T \int_{-1}^{+\infty} \frac{x^2}{1+x}d\mu = \infty \Big\}= \Big\{ 1 + \frac{1}{4} \int_0^T \int_{-1}^{+\infty} \frac{x^2}{1+x}d\mu = \infty \Big\} = $$

$$ \Big\{ \ln \Big( 1 + \frac{1}{4} \int_0^T \int_{-1}^{+\infty} \frac{x^2}{1+x}d\mu \Big) = \infty \Big\} \subset $$

$$ \Big\{ \int_0^T \int_{-1}^{+\infty} \ln \Big( 1+\frac{1}{4}\cdot \frac{x^2}{1+x} \Big)d\mu = \infty  \Big\} \subset  \{\mathcal{E}_{T}(M^d)=0 \}. $$
So we proved that $\{ \int_0^T \int_{-1}^{+\infty} \frac{x^2}{1+x}d\mu = \infty \} \subset \{\mathcal{E}_{T}(M^d)=0 \}$, for any local martingale $M$. It is clear that from this we can deduce as a particular case 
\begin{equation}
\Big\{ \int_0^T \int_{-1}^1 \frac{x^2}{1+x} d\mu = \infty \Big\} \subset \big\{ \mathcal E_{T}(M^1)=0 \big\}.
\end{equation}

Now it is time to prove the reverse set inclusion: 
$\{\mathcal{E}_{T}(M^1)=0 \} \subset \big\{ \int_0^T \int_{-1}^1 \frac{x^2}{1+x} d\mu = \infty \big\}.$ 

$$\mathcal{E}_{T}(M^1)\mathcal{E}^2_{T}\big(-\frac{1}{2} M^1\big)=\exp \Big\{
M_{T}^1 + \int_0^T \int_{-1}^1 \big[ \ln (1+x)-x \big] d\mu - M^1_T +$$
$$ \int_0^T \int_{-1}^1 \big[ 2\ln \big(1-\frac{x}{2}\big)+x \big]d\mu  \Big\}= \exp \Big\{ \int_0^T \int_{-1}^1 \ln \big[ (1+x)\big(1-\frac{x}{2}\big)^2 \big]d\mu \Big\}.$$
From the last equality and the supermartingale property of $\mathcal{E}( -\frac{1}{2}M^1) $ we deduce that

$$
\big\{ \mathcal{E}_{T }(M^1) = 0 \big\} \subset \Big\{ -\int_0^T \int_{-1}^1 \ln \big[ (1+x)\big(1-\frac{x}{2}\big)^2 \big]d\mu =\infty \Big\}.
$$
\\
Using Lemma 1 from Appendix we obtain $-\ln (1+x)\big( 1-\frac{x}{2}\big)^2 \leq \frac{2x^2}{1+x}$ and this gives us an inclusion:

$$\{ \mathcal{E}_{T}(M^1)=0\}\subset \Big\{ \int_0^T \int_{-1}^1 \frac{2x^2}{1+x}d\mu  = \infty \Big\}=\Big\{ \int_0^T \int_{-1}^1 \frac{x^2}{1+x}d\mu = \infty \Big\}$$
which with $(3)$ implies the equality $(1)$.

\
\\
Now we prove the set equality $\{ \mathcal{E}_{T}(M^2)=0\} = \big\{ \int_0^T \int_{(1;+\infty)} \frac{x^2}{1+x} d\nu = \infty \big\}$.
\\
It follows from Jacod \cite{1}, that 
$\big\{ \int^T_0 \int_1^{+\infty} \frac{x^2}{1+|x|}d\nu =\infty \big\}\subset \big\{ \mathcal{E}_T (M^2)=0 \big\}$. But it is clear that $\int^T_0 \int_1^{+\infty} \frac{x^2}{1+|x|}d\nu = \int^T_0 \int_1^{+\infty} \frac{x^2}{1+x}d\nu $, because $x\geq 1$. So we have $\big\{ \int^T_0 \int_1^{+\infty} \frac{x^2}{1+x}d\nu =\infty \big\}\subset \big\{ \mathcal{E}_T (M^2)=0 \big\}$. For the reverse inclusion it is clear that $\mathcal{E}_{T}(M^2)=\exp \big\{ -\int_0^{T}\int_{1}^{+\infty}x d\nu + \int_0^{T}\int_{1}^{+\infty} \ln (1+x) d\mu \big\}$ and from this we deduce $\{ \mathcal{E}_{T}(M^2)=0\}\subset \big\{ \int_0^{T}\int_{1}^{+\infty}x d\nu = \infty \big\}$. For $x\geq 1$ the inequality $x\leq \frac{2x^2}{1+x}$ holds true, which implies that $\big\{ \int_0^{T}\int_{1}^{+\infty}x d\nu = \infty \big\} \subset \big\{ \int_0^{T}\int_{1}^{+\infty}\frac{x^2}{1+x} d\nu = \infty \big\}$. So we will have inclusion
$
\{ \mathcal{E}_{T}(M^2)=0\} \subset \Big\{ \int_0^{T}\int_{1}^{+\infty}\frac{x^2}{1+x} d\nu = \infty \Big\}
$
and finally we get the set equality $(2)$. So the proof of part $(i)$ is completed.

\
\\
Now we shall prove part $(ii)$ and part $(iii)$ of Theorem 1. To prove part $(ii)$ we need the set equality 
\begin{equation}
\Big\{ \int^T_0 \int_{-1}^{1}\frac{x^2}{1+x}d\mu =\infty \Big\} = \Big\{ \int^T_0 \int_{-1}^{1}\frac{x^2}{1+x}d\nu =\infty \Big\}
\end{equation}
and for part $(iii)$ 
\begin{equation}
\Big\{ \int^T_0 \int_1^{+\infty}\frac{x^2}{1+x}d\mu =\infty \Big\} = \Big\{ \int^T_0 \int_1^{+\infty}\frac{x^2}{1+x}d\nu =\infty \Big\}.
\end{equation}
Inequality $\frac{(\bigtriangleup M_{\sigma})^2}{1+\bigtriangleup M_{\sigma}}1_{\{|\bigtriangleup M_{\sigma}|\leq 1\} }\leq \frac{1}{1+\bigtriangleup M_{\sigma}}1_{\{|\bigtriangleup M_{\sigma}|\leq 1\} }$ and the integrability condition $E\frac{1}{1+\bigtriangleup M_{\sigma}}1_{\{|\bigtriangleup M_{\sigma}|\leq 1\} }<\infty$ from part $(ii)$ gives us possibility to use Theorem 2.6.1 from \cite{4} to obtain $(4)$.
\\
By the same manner for $(iii)$ if we use inequality $\frac{x^2}{1+x}\leq x$ for $x\geq 1$, condition $E\bigtriangleup M_{\sigma}<\infty$ from part $(iii)$ and Theorem 2.6.1 from \cite{4}, we obtain $(5)$.
\qed

\
\\
{\bf 4. Appendix.}

\
\\
{\bf Lemma 1.} $k(x)=\frac{2x^2}{1+x}+\ln(1+x)(1-\frac{1}{2}x)^2\geq 0$ for any $x\in (-1;1)$.
\
\\
{\it Proof.} 
$$k'(x) = \frac{4x(1+x)-2x^2}{(1+x)^2}+\frac{(1-\frac{1}{2}x)^2-(1+x)(1-\frac{1}{2}x)}{(1+x)(1-\frac{1}{2}x)^2}=$$
$$\frac{2x^2 + 4x}{(1+x)^2} - \frac{3x}{(1+x)(2-x)}=\frac{-2x^3-3x^2+5x}{(1+x)^2(2-x)}=
\frac{x(2x+5)(1-x)}{(1+x)^2(2-x)}.$$
It is obvious that $k'(0)=0$, $k'(x)<0$ when $x\in (-1;0)$ and $k'(x)>0$ when $x\in (0;1)$. So $x=0$ is a minimum point and because $k(0)=0$, we can deduce that $k(x)\geq 0$ for $x\in (-1;1)$.
\qed


\newpage




\begin{thebibliography}{99}

\bibitem{1}
J.Jacod. Calcul Stochastique et Problemes de Martingales, 
Vol. 714 of {\it Lecture Notes in Mathematics}, Springer-Verlag, Berlin Heidelberg New York, 1979.

\bibitem{2}
N. Kazamaki. { \it Continuous Exponential Martingales and BMO }, Vol. 1579 of
                { \it Lecture Notes in Mathematics }, Springer, Berlin-Heidelberg, 1994.

\bibitem{3}
M. Larsson, J. Ruf.  Stochastic Exponentials and Logarithms on Stochastic Intervals - A Survey*, { \it Journal of Mathematical Analysis and Applications} Vol. 476, Issue 1, Issue on Stochastic Differential Equations, Stochastic Algorithms, and Applications, (2019)

\bibitem{4}
R. Sh. Liptser, A. N. Shiryaev. Theory of Martingales, 1986.


\end{thebibliography}
\end{document}